\documentclass[11pt,oneside,pdftex]{amsart}

\usepackage{geometry}   
\usepackage{wrapfig}
\usepackage{color,amsmath,amsfonts,amssymb,amsthm,color,mathtools}
\usepackage{cancel}
\usepackage{tensor}
\usepackage{tikz-cd}
\usepackage{bbm}
\usepackage{mathrsfs}
\usepackage{subfiles}
\usepackage{fancyhdr}
\usepackage{graphicx}
\usepackage{enumitem}
\usepackage{hyperref}

\newtheorem{thm}{Theorem}[section]
\newtheorem{lemma}[thm]{Lemma}
\newtheorem*{lemma*}{Lemma}
\newtheorem{prop}[thm]{Proposition}
\newtheorem{crol}[thm]{Corollary}

\newtheorem{conj}[thm]{Conjecture}

\theoremstyle{definition}
\newtheorem{defn}[thm]{Definition}
\newtheorem{constr}[thm]{Construction}

\newtheorem{rmk}[thm]{Remark}

\newtheorem{conv}[thm]{Convention}

\newcommand{\ot}{\widetilde}

\newcommand{\inv}{^{-1}}

\DeclareMathOperator{\im}{im}
\DeclareMathOperator{\rank}{rank}

\DeclareMathOperator{\interior}{int}
\DeclarePairedDelimiter{\abs}{\lvert}{\rvert}

\DeclarePairedDelimiterX{\inner}[2]{\langle}{\rangle}{#1 , #2}
\DeclarePairedDelimiterX{\pres}[2]{\langle}{\rangle}{#1 \delimsize\vert #2}

\newcommand{\zetminus}{ - }
\newcommand{\homology}{\ensuremath{{\sf{H}}}}

\begin{document}
%\lhead{left header}
%\rhead{right header}

\title{Nonpositive Towers in Bing's Neighbourhood }
\author[Max Chemtov]{Max Chemtov}
          \address{Dept. of Math. \& Stats.\\
                    McGill Univ. \\
                    Montreal, QC, Canada H3A 0B9 }
          \email{max.chemtov@mail.mcgill.ca}

\author[D.~T.~Wise]{Daniel T. Wise}
          \address{Dept. of Math. \& Stats.\\
                    McGill Univ. \\
                    Montreal, QC, Canada H3A 0B9 }
          \email{wise@math.mcgill.ca}
%{20F67, 20F65, 20E06}
\date{\today}
\thanks{Research supported by NSERC}

\begin{abstract}Every 2-dimensional spine of an aspherical 3-manifold has the nonpositive towers property, but every  collapsed 2-dimensional spine of a 3-ball containing a 2-cell has an immersed sphere.
\end{abstract}

\maketitle

%this is nuts
%meh, it keeps things organized
%but you cannot search the tex file for stuff..
%fair enough. i've only used tex for single-section things or solo/group assignments, but for one long document, that makes sense

\section{Introduction}

\begin{defn}
    A  2-complex $X$ has \emph{nonpositive immersions}
    if for every combinatorial immersion $Y\rightarrow X$ with $Y$ compact and connected, either $\chi(Y)\leq 0$ or 
    $Y$ is contractible.
    
    A  2-complex $X$ has \emph{nonpositive towers}
    if for every tower map $Y\rightarrow X$ with $Y$ compact and connected, either $\chi(Y)\leq 0$ or 
    $Y$ is contractible.
\end{defn}
There are many variations:
For instance, one can generalize to combinatorial near-immersions, or relax to $\pi_1Y=1$ or  $\chi(Y)\leq 1$,
and there also variations requiring  $\chi(Y)\leq -c|Y|$ for some ``size'' $|Y|$ of $Y$.
Note that nonpositive immersions implies nonpositive towers.
The main consequences of nonpositive immersions hold for nonpositive towers. E.g.\ if $X$ has nonpositive towers then $\pi_1X$ is locally indicable.
These ideas have promise as a contextualizing framework towards Whitehead's asphericity conjecture, as well as towards understanding  coherence.

In 
\cite{WiseNonPositiveCoherence} 
it was shown % we showed
that every aspherical 3-manifold with nonempty boundary
has a spine  with nonpositive immersions.
This utilized that there exists a spine with no near-immersion of a 2-sphere \cite{CorsonTrace2000}.

In this note, we observe the following failure, which is a special case of Proposition~\ref{prop:sphere_in_bing_neighbour}:% of nonpositive immersions:
\begin{thm}\label{thm:generalized bing}
Every  collapsed spine of a simply-connected 3-manifold containing a disc has an immersed sphere.
\end{thm}

\begin{figure}
    \centering
  \includegraphics[width=14cm]{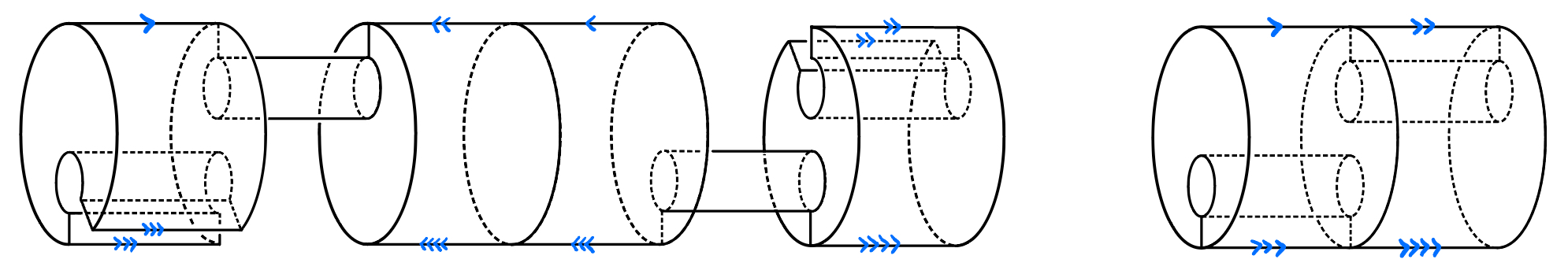}
    \caption{ An immersed sphere (left) in Bing's house (right).}
    \label{fig:bings_sphere}
\end{figure}

Bing's ``house with two rooms'' provides such a spine. It is obtained from a 3-ball divided into two rooms by a pair of collapses, corresponding to entering the left room from the right side of the house and entering the right room from the left side. See Figure~\ref{fig:bings_sphere}.

W.~Fisher also found examples of the failure of nonpositive immersions in other contractible 2-complexes: the  Miller-Schupp balanced presentations of the trivial group \cite{FisherWilliam2022}. There are thus two sources of counter-examples to the conjecture that contractible 2-complexes have nonpositive immersions \cite[Conj~1.7]{WiseNonPositiveCoherence}.

In \cite{FreedmanNguyen-Phan2021} it is shown that for $n\geq 3$, every PL $n$-manifold $M$ with $\partial M\neq \emptyset$ has a spine $X$ such that $\partial M\rightarrow X$ is an immersion, and moreover, such spines are generic among all spines.
So Theorem~\ref{thm:generalized bing} is a  variant of the simplest instance of their result. However it is a counterpoint to the following statement (see Theorem~\ref{nonpositive towers for aspherical 3-manifold spines}) which is our main motive:

\begin{thm}
Every 2-dimensional spine of an aspherical 3-manifold has nonpositive towers.
\end{thm}

Thus nonpositive immersions does not always hold for a natural family of contractible complexes which nevertheless have nonpositive towers, so we are motivated to refocus on:
\begin{conj}
Every contractible 2-complex has nonpositive towers.
\end{conj}
For instance:
\begin{prop} A contractible 2-complex with two 2-cells has nonpositive towers.
\end{prop}
\begin{proof}
As a contractible 2-complex $X$ admits no nontrivial connected covering map, any tower map $Y\rightarrow X$ must begin with a subcomplex $X'\subset X$, so $X'$ has at most one 2-cell.
But  $X'$ has nonpositive immersions (hence nonpositive towers) by \cite{LouderWilton2014,HelferWise2015}.
\end{proof}

{\bf Acknowledgement:} We are extremely grateful to Grigori Avramidi for helpful comments and the valuable reference \cite{FreedmanNguyen-Phan2021}. We tremendously appreciate the referee's useful feedback in improving the exposition.

\section {Thickenings and nonpositive towers}

%((hide after adjuting: we just assume towers use finite subcomplexes, since the top of the tower compact... so no need for towers... also we remember to add howie reference I'll send bibtex file now))

%In this section, we construct a thickening of a 2-complex, and use it to find a broad class of 2-complexes with nonpositive towers.

\begin{defn}
    Let $X$ be a connected 2-complex. A \emph{tower map} $Y\to X$ is a finite composition of covering maps and subcomplex embeddings, such that $Y$ and the domain of each embedding are compact connected 2-complexes.
\end{defn}

Tower maps arose in Papakyriakopolous' classical 3-manifold proofs, and also arose  naturally in 
one-relator group theory 
\cite{Howie87}.

\begin{defn}
    A 2-complex $X$ has \emph{nonpositive towers} if, for any tower map $Y\to X$, either $\chi(Y)\leq 0$ or $Y$ is contractible.
\end{defn}

The goal of this section is to prove the nonpositive tower property for an aspherical 2-complex $X$ embedded in a 3-manifold $M$. The idea of the proof is to consider a manifold ``thickening'' $T$ of $X$ in $M$, which deformation retracts to $X$. The asphericity of $X$ ensures that $\partial T$ has no 2-sphere, which in turn ensures $\chi(X)\leq 0$. Asphericity is preserved by towers, since it is preserved by both covering maps and subcomplexes. The latter is ensured by a simple argument using the Sphere Theorem.

\begin{conv}
    We work in the category of PL-manifolds: all submanifolds, as well as cells of complexes embedded in a manifold, are assumed to be PL-embedded. This avoids pathologies like Alexander's horned sphere.
\end{conv}

\begin{constr}
Let $X$ be a locally finite 2-complex. A \emph{thickening} of $X$ is a 3-manifold $T=T(X)$ with boundary, and a continuous map $\Theta:T\to X$, constructed as follows:
\begin{itemize}
    \item Let $T^0$ be a disjoint union of closed 3-balls, one for each vertex in $X$, and let $\Theta$ map each ball to its corresponding vertex.
    
    \item For each edge $e$ in $X$, define $T(e)\cong [0,1]\times D^2$, and identify  $\{0\}\times D^2$ and  $\{1\}\times D^2$ with discs on the boundary of the components of $T^0$ corresponding to the endpoints of $e$. Let $\Theta$ map $(0,1)\times D^2$ onto $\interior(e)$. The resulting complex is $T^1$.
    
    We require that each $T(e)$ embeds and that $T(e_1)\cap T(e_2)=\varnothing$ for $e_1\neq e_2$.
    
    \item For each disc $F$ in $X$, define $T(F)\cong D^2\times [0,1]$, and identify the outer cylinder $S^1\times [0,1]$ with an embedded cylinder on the boundary of $T^1$ which is mapped by $\Theta$ to the attaching loop of $F$ in $X$. Let $\Theta$ map $\interior(D^2)\times [0,1]$ onto $\interior(F)$. The resulting complex is $T$.
    
    We require that each $T(F)$ embeds and that $T(F_1)\cap T(F_2)=\varnothing$ for $F_1\neq F_2$.
    
\end{itemize}
For a subcomplex $A\subseteq X$, we use $T(A)\subseteq T(X)$ to denote the thickening of $A$ to $\Theta\inv(A)$ induced by the thickening of $X$.

\end{constr}

\begin{rmk}
    If a thickening $T$ of $X$ exists, then by construction, $T$ is a 3-manifold with boundary. Furthermore, there is a PL-embedding $X\hookrightarrow \interior(T(X))$ such that $T$ deformation retracts to $X$, with $T(A)$ retracting to $A$ for any subcomplex $A\subseteq X$.
\end{rmk}
\begin{rmk}
    If $X$ PL-embeds in a 3-manifold $M$, then we can take $T$ to be a small closed neighbourhood of $X$ in $M$. Let $\Theta:T(X)\to X$ be a retraction homotopic to the identity map $T(X)\to T(X)$. Then $T$ and $\Theta$ give a thickening of $X$.
\end{rmk}

%%%%%%%%%%%% FIX DEFINITION OF T(X), Theta

\begin{lemma} \label{deleting 2-cells preserves asphericity}
    Let $X$ be a locally finite aspherical connected 2-complex that PL-embeds in a 3-manifold. Then removing a collection of 2-cells $\{D_i\}$ in $X$ results in a new 2-complex $Y$ which is also aspherical.
\end{lemma}

\begin{proof}
    If $T(X)$ is non-orientable, we can consider an orientable double-cover $\widehat{T(X)}$. This induces a double-cover $\widehat X$ of $X$, which is locally finite, aspherical, connected, and PL-embeds in $T(\widehat X) = \widehat {T(X)}$. Consider the induced orientable double-cover $\widehat Y \subseteq \widehat X$ of $Y$, and note that $\widehat Y$ can be obtained from $\widehat X$ by deleting a collection of 2-cells. It suffices to prove that $\widehat Y$ is aspherical, since this would imply that $Y$ is aspherical. Since the non-orientable case with $X$ and $Y$ reduces to the orientable case with $\widehat X$ and $\widehat Y$, we can assume without loss of generality that $T(X)$ is orientable.

    Let $Y$ be a 2-complex, and let $X = Y\cup\left(\bigcup\limits_i D_i\right)$, where each 2-cell $D_i$ is attached to $Y$ along $\partial D_i$. We assume $X$ is a locally finite aspherical connected 2-complex that PL-embeds in an orientable manifold, and prove that $Y$ is aspherical.
    
    Suppose for contradiction that $Y$ is not aspherical. Let $(T(X), \Theta)$ be an orientable thickening of $X$. Then $\pi_2(\interior(T(Y)))\neq 0$. Since $\interior(T(Y))$ is orientable, $\interior(T(Y))$ has an embedded essential 2-sphere $S$
    by %Papakyriakopoulos'
    the Sphere Theorem \cite[Thm~3.8]{Hatcher3Man}. Since $X$ is aspherical, $S$ bounds a contractible submanifold $B\subseteq \interior(T(X))$
     by \cite[Prop~3.10]{Hatcher3Man}. Note that $T(X)\zetminus S$ has two connected components: $\interior(B)$ and the component $C$ containing $\partial T(X)$.
    
    For any $D_i$, we have $\Theta\inv(\interior(D_i))\cap S=\varnothing$, since $S\subseteq T(Y)$ and $\interior(D_i)\cap Y=\varnothing$. Since $\Theta\inv(\interior(D_i))$ is connected, it must lie either entirely in $\interior(B)$ or entirely in $C$. By construction of  $T(X)$, we know that $\Theta\inv(\interior(D_i))\cap \partial T(X)\neq\varnothing$. So $\Theta\inv(\interior(D_i))  \subset C$.
    
    Since this is true for all $D_i$, we have $\interior(B)\subseteq T(X)\zetminus\bigcup\limits_i \Theta\inv(\interior(D_i)) = T(Y)$. Since $\interior(B)$ is an open submanifold of $T(Y)$, it is contained in $\interior(T(Y))$. Therefore, $S$ bounds a contractible submanifold $B$ of $\interior(T(Y))$, contradicting that $S$ is essential in $\interior(T(Y))$.
\end{proof}

\begin{crol} \label{subcomplex preserves asphericity}
    Let $X$ be a locally finite aspherical connected 2-complex that PL-embeds in a 3-manifold. Then every subcomplex $Y\subseteq X$ is aspherical.
\end{crol}
\begin{proof}
    By Lemma~\ref{deleting 2-cells preserves asphericity}, removing 2-cells from $X$ yields another aspherical 2-complex. Since removing 0- and 1-cells also preserves asphericity,  every subcomplex of $X$ is aspherical.
\end{proof}

%\begin{lemma} \label{3-manifold euler characteristic}
%    Let $M$ be a compact connected orientable 3-manifold, with $c>=1$ components in $\partial M$.
%    Then $$1-b_1(M)+(c-1) \leq 1-b_1(M)+b_2(M) = \chi(M)=\frac{1}{2}\chi(\partial M)$$ Thus if $\partial M$ does not contain a 2-sphere then $\chi(M)\leq 0$. And if $b_1(M)=0$ then $\partial M = S^2$ \end{lemma}
%\begin{proof}
%    Let $\ot M$ be the manifold obtained by gluing two copies of $M$  along $\partial M$. Since $\ot M$ is a closed 3-manifold,  $ 2\chi(M) - \chi(\partial M) = \chi(\ot M) = 0$. So $\chi(M)=\frac{1}{2}\chi(\partial M)$.
%    This establishes the final equality, and the initial inequality holds  as $M$ is connected, so    any $c-1$ of the boundary 2-cycles are independent in $\homlogy_2(M)$. 
%    
%    The first claim holds since each component of $\partial M$ is an orientable surface. So     if  $\partial M$ contains no 2-sphere, every component of $\partial M$ has nonpositive $\chi$. So $\chi(M)=\frac{1}{2}\chi(\partial M) \leq 0$.
    %The second claim holds since  OOPS DIDNT WORK
%\end{proof}

\begin{lemma} \label{3-manifold has nonpositive euler characteristic}
    Let $M$ be a compact orientable 3-manifold with boundary. Then $\chi(M)=\frac{1}{2}\chi(\partial M)$. And if $\partial M$ does not contain a 2-sphere then $\chi(M)\leq 0$.
\end{lemma}
\begin{proof}
    Let $\ot M$ be the manifold obtained by gluing two copies of $M$  along $\partial M$. Since $\ot M$ is a closed 3-manifold,  $ 2\chi(M) - \chi(\partial M) = \chi(\ot M) = 0$. So $\chi(M)=\frac{1}{2}\chi(\partial M)$. Since $M$ is orientable, each component of $\partial M$ is an orientable surface. Since $\partial M$ contains no 2-sphere, every component of $\partial M$ has nonpositive $\chi$. So $\chi(M)=\frac{1}{2}\chi(\partial M) \leq 0$.
\end{proof}

\begin{lemma} \label{finite covers preserve nonpositive towers}
    Let $X'\to X$ be a finite-sheeted cover of a compact connected 2-complex X. Then $X$ has nonpositive towers if and only if $X'$ has nonpositive towers.
\end{lemma}
\begin{proof}
    Suppose  $X$ has nonpositive towers. Let $Y\to X'$ be a tower map. Then $Y\to X'\twoheadrightarrow X$ is a tower map, so either $\chi(Y)\leq 0$ or $Y$ is contractible. 
    
    Suppose that $X'$ has nonpositive towers. Let $t:Y\to X$ be a tower map. Let $n$ be the degree of the cover $p:X'\to X$. Then there is an induced tower map $Y'\to X'$, where $Y'$ is an $n$-sheeted cover of $Y$, and maps to $p\inv(t(Y))$. Either $\chi(Y')\leq 0$ or $Y'$ is contractible.
    
    If $\chi(Y')\leq 0$, then $\chi(Y)=\frac{1}{n}\chi(Y')\leq 0$. If $Y'$ is contractible, then $Y'$ is the universal cover of $Y$,  so $Y$ is a $K(\pi_1(Y), 1)$ complex with $\abs{\pi_1 Y}=n$. A nontrivial finite group does not have a compact $K(\pi,1)$, so $\pi_1 Y =1 $.
     Thus $Y=Y'$ is contractible.
    %    Since $Y$ was arbitrary, $X$ has nonpositive towers.
\end{proof}

\begin{thm} \label{nonpositive towers for aspherical 3-manifold spines}
    Let $X$ be an aspherical compact connected 2-complex that PL-embeds in a 3-manifold $M$. Then $X$ has nonpositive towers.
\end{thm}

\begin{proof}
    Since $X$ is compact, it PL-embeds in a thickening $T(X)$ in $M$ that is a compact sub-\hbox{3-manifold} with boundary. So without loss of generality, we can assume that $M=T(X)$ is a compact 3-manifold with boundary that deformation retracts to $X$.
    
    If $M$ is non-orientable,  we  consider an orientable double cover of $M$ that deformation retracts to a double cover $X'$ of $X$. By Lemma~\ref{finite covers preserve nonpositive towers}, it suffices to show that $X'$ has nonpositive towers. So without loss of generality, we can assume that $M$ is orientable.

    $X$ itself is either contractible or has $\chi(X)\leq 0$. Indeed, if $\partial M$ includes a 2-sphere $S$, then asphericity of $M$ ensures that $S$ bounds a contractible submanifold of $M$ \cite[Prop~3.10]{Hatcher3Man}. Since $M$ is connected, this submanifold must be $M$ itself, and so $M$ and $X$ are contractible. If $\partial M$ does not include any 2-spheres, then $\chi(X) = \chi(M) \leq 0$ by Lemma~\ref{3-manifold has nonpositive euler characteristic}.
    
    Note that any covering map $\widehat X \to X$ (with $\widehat X$ connected) extends to a covering map $\widehat M \to M$, where $\widehat X$ PL-embeds in $\interior(\widehat M)$. Then $\widehat X$ is locally finite, aspherical, connected, and PL-embeds in a 3-manifold.
    
    Let $X'$ be a compact connected subcomplex of $\widehat X$. By Corollary \ref{subcomplex preserves asphericity}, $X'$ is aspherical. And $X'$ also PL-embeds in $\interior(\widehat M)$. Therefore, $X'$ satisfies the same hypotheses as $X$, and is also either contractible or has $\chi(X')\leq 0$.
    
    Any tower map $Y\to X$ is a composition of maps $X'\hookrightarrow \widehat X \twoheadrightarrow X$, so we are done.
\end{proof}

\section{Bing's Neighbours}

In this section, we give an alternate construction of the thickening $T$ of a 2-complex $X$ embedded in a 3-manifold. The cell structure on $\partial T$ ``follows'' $X$ and provides an immersion $\partial T \to X$. When $X$ is simply-connected, collapsed, and contains a 2-cell, $\partial T$ is a union of 2-spheres. Thus examples like Bing's house fail to have nonpositive immersions.

\begin{defn}
    A 2-complex is \emph{collapsed} if no cell has a free face - i.e.\ no 0-cell has degree~1, and no 
    1-cell is incident to a single side of a 2-cell.
%   A \emph{free face} is an $i$-cell that appears exactly once, and in only one way among all the attaching maps of $(i+1)$-cells.
\end{defn}

\begin{constr}
Let $X$ be a compact connected collapsed 2-complex with no isolated vertex or edge that PL-embeds in a 3-manifold $M$. We give a construction for the thickening $T=T(X)$, yielding an explicit cell structure on $\partial T$ related to the cell structure on $X$:
\begin{enumerate}[leftmargin=.7cm]
    \item[$\partial T^0$] For each vertex $v$ in $X$, consider a small neighbourhood $N(v)\subseteq M$ of $v$. Add a 0-cell in each component of $N(v) - X$. The union of these 0-cells is $\partial T^0$.
    
    \item[$\partial T^1$] For each edge $e$ in $X$, consider a small neighbourhood $N(e)\subseteq M$ of $e$ containing the 0-cells in $\partial T^0$ associated with the endpoints of $e$. Add a 1-cell in each component of $N(e) - X$ that contains a 0-cell associated to each endpoint of $e$. That 1-cell joins those 0-cells. If $n$ sides of discs are incident with $e$ in $X$, then this process yields $n$ 1-cells parallel to $e$. The union of $\partial T^0$ with these 1-cells is $\partial T^1$.
    
    \item[$\partial T^2$] For each disc $d$ of $X$, consider a small neighbourhood $N(d) \subseteq M$ of $d$ containing the 1-cells in $\partial T^1$ associated with the $\partial d$. Add a 2-cell in each component of $N(d) - X$ containing 1-cells associated to all edges of $\partial d$. Attach this 2-cell to those 1-cells according to $\partial d$. This yields two 2-cells on opposite sides of $d$ in $M$. The union of $\partial T^1$ with these 2-cells is $\partial T$.
    
    %\item[$H^3$] Finally, attach a 3-cell filling the compact component of $M - H^2$ containing $X$. The resulting complex is $H$.
\end{enumerate}

\noindent The submanifold $T$ is the union of $\partial T$ and the component of $M\zetminus\partial T$ containing $X$.
\end{constr}

\begin{lemma} \label{properties of H}
    $T=T(X)$ deformation retracts to $X$. The retraction $r:T\to X$ induces an immersion $\partial T \to X$. If $X$ is simply-connected, then $\partial T$ is a union of 2-spheres.
\end{lemma}

\begin{figure}
    \centering
    \includegraphics[width=5cm]{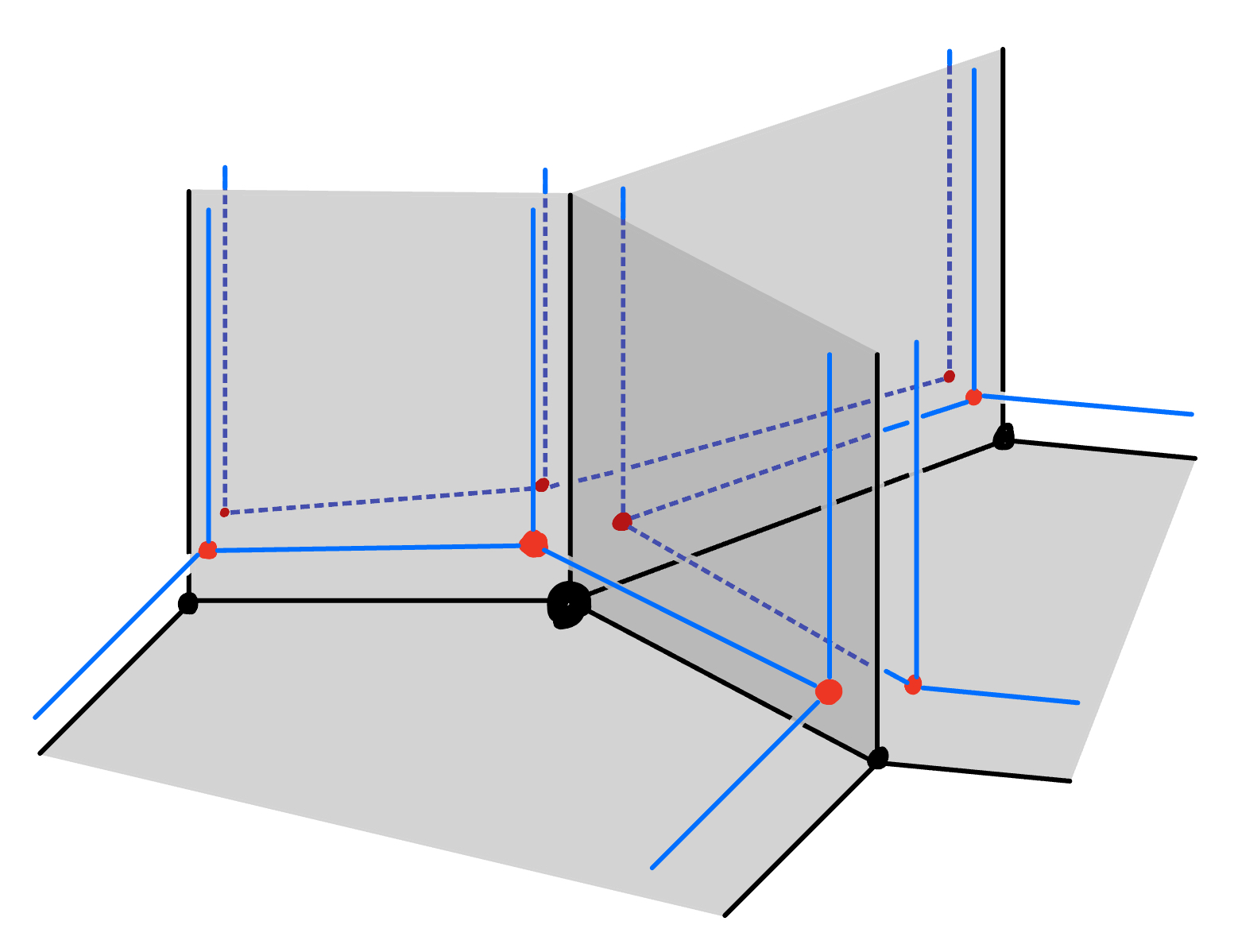}
    \caption{Part of $\partial T(X)$ for a complex $X$. The 0-cells of $\partial T$ are red, and the 1-cells of $\partial T$ are blue. The 2-cells of $\partial T$ run parallel to the grey discs in $X$, forming a ``bubble'' hovering around $X$.}
    \label{fig:bubble thickening}
\end{figure}

\begin{proof}
     Consider the map $\partial T\to X$ sending $i$-cells in $\partial T$ to their associated $i$-cells in $X$. \hbox{$T$ is homeomorphic} to the mapping cylinder of $\partial T\to X$, yielding a deformation retraction.
     
     %of $M$ that deformation retracts to $X$ (in fact, $T$ is the thickening $T$ defined in the previous section). One such deformation retraction $r$ is induced by mapping $i$-cells in $\partial T$ to their associated $i$-cells in $X$.
    
    Let $c_1$ and $c_2$ be closed $i$-cells in $\partial T$ with $c_1\cap c_2\neq\varnothing$. Suppose  $r(c_1)=r(c_2)=c$. Then $c_1\cup c_2$ is a connected subset of the neighbourhood $N(c)$ used in the construction of $\partial T$. At most one $i$-cell was added for each  component of $N(c)$, so  $c_1=c_2$. Thus $r|_{\partial T}:\partial T\to X$ is an immersion.
    
    The map $\partial T\to X$ can also be seen geometrically to be an immersion: in terms of Figure~\ref{fig:bubble thickening}, each cycle of 2-cells around a vertex in $\partial T$ is mapped to a cycle of discs in $X$ associated to a corner of $M-X$.
    
    %I think I prefer to understand the immersion by thinking geometrically: A cycle of polygons around a vertex in $H$ maps to a cycle of polygons around a 3-dimensional corner in $M-X$ and thus immerses in $X$.
    
    Suppose  $\pi_1X=1$. Then $T$ is a compact orientable 3-manifold with $\rank(\homology_1(T))=0$, so
    $\rank(\homology_1(\partial T))=0$ by ``half lives, half dies'' \cite[Lem~3.5]{Hatcher3Man}.
    Thus $\partial T$ is  a union of 2-spheres.
\end{proof}

\begin{lemma} \label{good subcomplex with no free edges}
    Let $X$ be a compact 2-complex with a collapsed subcomplex $Y$ containing a disc. Then $X$ has a connected collapsed subcomplex $Y'$ with no isolated vertex or edge, whose inclusion map $Y'\hookrightarrow X$ is $\pi_1$-injective.
\end{lemma}

\begin{proof}
    If $X$ contains a free $i$-face, delete that face and its attached $(i+1)$-cell. This deletion does not affect $\pi_1$. Repeat this process until the remaining subcomplex is collapsed. Delete all isolated edges leaving a disjoint union $X'$ of collapsed components $X'_i$  without isolated vertices or edges. Each  $X'_i\hookrightarrow X$ is $\pi_1$-injective. As  the disc from $Y$ must be contained in some $X'_i$, we let $Y'=X'_i$.
\end{proof}

% \begin{proof}
%     Let $Y\subseteq X$ have no free edges and be more than a single point. Suppose that $Y$ has an essential loop $\gamma$ given by a freely reduced word in its 1-cells, which is null-homotopic in $X$. Then there is a reduced disc diagram $D\to X$ with $\partial D$ mapping to $\gamma$. Since $D$ is reduced, the only free edges of its image in $X$ are along $\gamma$. Since the edges in $\gamma$ are each incident with at least two sides of 2-cells in $Y$, we get a new subcomplex $Y\cup\im(D)$ of $X$ with no free edges. Note that $\gamma$ is null-homotopic in $Y\cup\im(D)$.
    
%     Repeating this process a finite number of times (since $X$ is compact), we can continue adding disc diagrams to $Y$ until we get a complex $Y'$, where every essential loop is also essential in $X$. This means that $\ker(\pi_1(Y')\to\pi_1(X))=0$, so $Y'$ is our desired subcomplex.
% \end{proof}

\begin{prop}\label{prop:sphere_in_bing_neighbour}
    Let $X$ be a 2-spine of a simply-connected 3-manifold. Then $X$ has an immersed 2-sphere if and only if $X$ has a collapsed subcomplex containing a disc.
\end{prop}
\begin{proof}
    Suppose $S^2\looparrowright X$ is a combinatorial immersion. Then $\im(S^2\looparrowright X)$ is a collapsed subcomplex containing a disc.
    
    Conversely, suppose $X$ has a collapsed  subcomplex $Y$ containing a disc. By Lemma~\ref{good subcomplex with no free edges}, we can assume $Y$ has no isolated vertex or edge, and $\pi_1Y\leq\pi_1X=1$. Then by Lemma \ref{properties of H}, $\partial T(X)$ is a union of 2-spheres that immerses in $X$.
\end{proof}

%\noindent {\bf Acknowledgement:} We are grateful to the referee for many helpful  corrections.}

%this is the usual way... not for overleaf
%\bibliographystyle{alpha}
%%%\bibliography{C:/Users/Dani/Dropbox/papers/wise}
%\bibliography{Sections/wise}

\bibliographystyle{plain} % We choose the "plain" reference style
\bibliography{wise} % Entries are in the refs.bib file

%\medskip
%\printbibliography
\end{document}